\documentclass{article}
\usepackage[T2A]{fontenc}
\usepackage[cp1251]{inputenc}
\usepackage[english,russian]{babel}
\usepackage[tbtags]{amsmath}
\usepackage{amsfonts,amssymb,mathrsfs,amscd}
\numberwithin{equation}{section}
\newtheorem{remark}{Remark}
\newtheorem{lemma}{Lemma}
\newtheorem{theorem}{Theorem}
\newtheorem{definition}{Definition}

\newcommand{\sign}{\mathop{\rm sign}}
\newcommand{\res}{\mathop{\rm res}}
\renewcommand{\span}{\mathop{\rm span}}

\renewcommand{\d}{\mathop{\rm def}}

\begin{document}
\begin{Large}
\thispagestyle{empty}
\begin{center}
{\bf Direct and inverse problems for a periodic problem with non-local potential\\
\vspace{5mm}
V. A. Zolotarev}\\
\end{center}
\vspace{5mm}

{\small {\bf Abstract.} Sum of a second derivative operator  with periodic boundary conditions and an integral operator of rank one (non-local potential) is studied in this manuscript. Not only spectral analysis is conducted for this operator but the inverse problem is solved, viz., a method of the reconstruction of non-local potential by the spectral data is presented. Class of spectral data that give the solution to the inverse problem is described.}
\vspace{5mm}


{\it Key words}: non-local potential, direct problem, inverse problem, periodic boundary problem.
\vspace{5mm}

Spectral analysis of the Sturm -- Liouville operator on a finite interval is well-studied and has multiple applications \cite{1} -- \cite{3}. Operators with non-local potential emerged for the first time in the study of behavior of particles on a crystalline surface \cite{4}. Such operators are sums of a second derivative operator and an integral operator. In the case when the differential operator has simple spectrum, direct and inverse spectral problems are solved in the papers \cite{5}, \cite{6}.

This manuscript deals with the direct and inverse problems for the operator
$$(Ly)(x)=-y''(x)+\alpha\int\limits_0^\pi y(t)v(t)dtv(x),\quad x\in(0,\pi),$$
where $\alpha\in\mathbb{R}$, $v\in L^2(0,\pi)$, with the domain $y\in W_2^2(0,\pi)$ and
$$y(0)=y(\pi),\quad y'(0)=y'(\pi).$$
Initial second derivative operator in $L^2(0,\pi)$ with the domain defined by these boundary conditions has discrete spectrum of multiplicity two. In the manuscript description of the spectrum of operator $L$ is presented and it is shown that spectral multiplicity of the eigenvalues of operator $L$ can be equal to 1, 2, 3. It is shown that the number $\alpha$ and the real function $v(x)$ are unambiguously defined by the three spectra.
\vspace{5mm}

\section{Characteristic function. Eigenfunctions}\label{s1}

{\bf 1.1.} Let us denote by $L_0$ a self-adjoint non-negative operator in $L^2(0,\pi)$,
\begin{equation}
\left(L_0y\right)(x)\stackrel{\d}{=}-\frac{d^2}{dx^2}y(x),\label{eq1.1}
\end{equation}
domain of which is defined by the boundary conditions
\begin{equation}
\mathfrak{D}\left(L_0\right)=\left\{y\in W_2^2(0,\pi):y(0)=y(\pi),y'(0)=y'(\pi)\right\}.\label{eq1.2}
\end{equation}
Eigenfunctions of the operator $L_0$ satisfy the equation
$$-y''(\lambda,x)=zy(\lambda,x)\quad\left(z=\lambda^2\right),$$
general solution to which is given by
$$y(\lambda,x)=A\cos\lambda x+B\frac{\sin\lambda x}\lambda.$$
Boundary conditions for $y(\lambda,x)$ imply the system of equations
$$\left\{
\begin{array}{lll}
{\displaystyle A(\cos\lambda\pi-1)+B\frac{\sin\lambda\pi}\lambda=0;}\\
-A\lambda\sin\lambda\pi+B(\cos\lambda\pi-1)=0,
\end{array}\right.$$
for $A$ and $B$, solutions to which are non-trivial only when its determinant $\Delta(0,\lambda)$ vanishes, $\Delta(0,\lambda)=0$,
$$\Delta(0,\lambda)=(\cos\lambda\pi-1)^2+\sin^2\lambda\pi=2(1-\cos\lambda\pi)=0.$$
The function
\begin{equation}
\Delta(0,\lambda)\stackrel{\d}{=}2(1-\cos\lambda\pi)=2-e^{i\lambda\pi}-e^{-i\lambda\pi}\label{eq1.3}
\end{equation}
is said to be the characteristic function of the operator $L_0$ \eqref{eq1.1}, \eqref{eq1.2}. It has the properties
\begin{equation}
\Delta(0,-\lambda)=\Delta(0,\lambda);\quad\Delta^*(\lambda)=\Delta(\lambda)\label{eq1.4}
\end{equation}
where
\begin{equation}
f^*\stackrel{\d}{=}\overline{f\left(\bar\lambda\right)}.\label{eq1.5}
\end{equation}
Zeroes $\lambda_k(0)=2k$ ($k\in\mathbb{Z}_+$) of the function $\Delta(0,\lambda)$ have the algebraic multiplicity 2 ($\Delta\left(0,\lambda_k(0)\right)=\Delta'\left(0,\lambda_k(0)\right)=0$) and define spectrum of the operator $L_0$,
\begin{equation}
\sigma\left(L_0\right)=\left\{z_k(0)=\lambda_k^2(0)=4k^2:k\in\mathbb{Z}_+\right\}.\label{eq1.6}
\end{equation}
To the eigenvalue $z_0(0)=0$ ($k=0$) there corresponds the one-dimensional eigenspace, it is defined by the normalized eigenfunction ${\displaystyle u(0,x)=\frac1{\sqrt\pi}}$. To the eigenvalue $z_k(0)=4k^2$ (to the root $\lambda_k(0)=2k$, $k\in\mathbb{N}$) there corresponds the two-dimensional eigenspace
\begin{equation}
G_k=\span\left\{\sqrt\frac2\pi\cos 2kx,\sqrt\frac2\pi\sin 2kx\right\}.\label{eq1.7}
\end{equation}
Thus all the eigenvalues $z_k(0)$ ($k\in\mathbb{N}$) have spectral multiplicity 2, except for the number $z_0(0)$ which has spectral multiplicity 1. The totality
\begin{equation}
\frac1{\sqrt\pi},\quad\left\{\sqrt{\frac2\pi}\cos 2kx,\sqrt{\frac2\pi}\sin 2kx\right\}_{k\in\mathbb{N}}\label{eq1.8}
\end{equation}
is an orthonormal basis in $L^2(0,\pi)$.
\vspace{5mm}

{\bf 1.2.} Consider the self-adjoint operator
\begin{equation}
L=L(\alpha,v)\stackrel{\d}{=}L_0+\alpha\langle\cdot,v\rangle v\label{eq1.9}
\end{equation}
which is a one-dimensional perturbation of the operator $L_0$ \eqref{eq1.1}, \eqref{eq1.2} where $\alpha\in\mathbb{R}$ and $v\in L^2(0,\pi)$. Domains of the operator $L$ and $L_0$ coincide, $\mathfrak{D}(L)=\mathfrak{D}\left(L_0\right)$ \eqref{eq1.2}.

\begin{remark}\label{r1.1}
On account of renormalization, we assume that either $\alpha=\pm1$ ($v\rightarrow\sqrt{|\alpha|}v$) or $\|v\|=1$ ($\alpha\rightarrow\alpha\|v\|^2$).
\end{remark}

Solving the equation
\begin{equation}
-y''(x)-\lambda^2y(x)=f(x)\label{eq1.10}
\end{equation}
by the method of variation of arbitrary constants, we obtain that the function
\begin{equation}\label{eq1.11}
\begin{array}{ccc}
\displaystyle{y(\lambda,x)=a(\lambda)\left\{\frac{\sin\lambda x}\lambda+\frac{\sin\lambda(\pi-x)}\lambda\right\}-\frac2{\Delta(0,\lambda)}\left\{\int\limits_0^x\frac{\sin\lambda(x-t)}\lambda f(t)dt\right.}\\
\displaystyle{\left.+\frac{\sin\lambda(\pi-x)}\lambda\int\limits_0^x\cos\lambda tf(t)dt+\cos\lambda x\int\limits_x^\pi\frac{\sin\lambda(\pi-t)}\lambda f(t)dt\right\}.}
\end{array}
\end{equation}
is the solution to equation \eqref{eq1.10} and satisfies the first boundary condition $y(\lambda,0)=y(\lambda,\pi)$. Hence it follows that the eigenfunction $u(\lambda,x)$ of the operator $L$ \eqref{eq1.9},
$$-u''(\lambda,x)+\alpha\langle u,v\rangle v(x)=\lambda^2u(\lambda,x),$$
is the solution to the integral equation
\begin{equation}
\begin{array}{ccc}
\displaystyle{u(\lambda,x)=a(\lambda)\left\{\frac{\sin\lambda x}\lambda+\frac{\sin\lambda(\pi-x)}\lambda\right\}+\frac{2b(\lambda)}{\Delta(0,\lambda)}\left\{\int\limits_0^x\frac{\sin\lambda(x-t)}\lambda v(t)dt\right.}\\
\displaystyle{\left.+\frac{\sin\lambda(\pi-x)}\lambda\int\limits_0^x\cos\lambda tv(t)dt+\cos\lambda x\int\limits_x^\pi\frac{\sin\lambda(\pi-t)}\lambda v(t)dt\right\}}\label{eq1.12}
\end{array}
\end{equation}
and satisfies the boundary condition $u(\lambda,0)=u(\lambda,\pi)$ where $b(\lambda)=\alpha\langle u,v\rangle$. Scalarly multiplying equality \eqref{eq1.12} by $\alpha v(x)$, we obtain
\begin{equation}
\alpha a(\lambda)W(\lambda)+\frac{2b(\lambda)}{\Delta(0,\lambda)}\left\{\alpha F(\lambda)+\alpha F^*(\lambda)+\alpha G(\lambda)-\frac{\Delta(0,\lambda)}2\right\}=0\label{eq1.13}
\end{equation}
where
$$W(\lambda)=\int\limits_0^\pi\frac{\sin\lambda x}\lambda\bar v(x)dx+\int\limits_0^\pi\frac{\sin\lambda(\pi-x)}\lambda\bar v(x)dx;$$
\begin{equation}
F(\lambda)=\int\limits_0^\pi\frac{\sin\lambda(\pi-x)}\lambda\int\limits_0^x\cos\lambda tv(t)dt\bar v(x)dx;\label{eq1.14}
\end{equation}
$$G(x)=\int\limits_0^\pi\int\limits_0^x\frac{\sin\lambda(x-t)}\lambda v(t)dt\bar v(x)dx.$$
Formula \eqref{eq1.12} yields
$$u'(\lambda,x)=a(\lambda)\{\cos\lambda x-\cos\lambda(\pi-x)\}+\frac{2b(\lambda)}{\Delta(0,\lambda)}\left\{\int\limits_0^x\cos\lambda(x-t)v(t)dt\right.$$
$$\left.-\cos\lambda(\pi-x)\int\limits_0^x\cos\lambda tv(t)dt-\sin\lambda x\int\limits_x^\pi\sin\lambda(\pi-t)v(t)dt\right\}$$
and according to the second boundary condition, $u'(\lambda,0)=u'(\lambda,\pi)$,
\begin{equation}
a(\lambda)\Delta(0,\lambda)+\frac{2b(\lambda)}{\Delta(0,\lambda)}N(\lambda)=0\label{eq1.15}
\end{equation}
where
\begin{equation}
N(\lambda)=\int\limits_0^\pi(\cos\lambda t-\cos\lambda(\pi-t))v(t)dt.\label{eq1.16}
\end{equation}
From the relations \eqref{eq1.13}, \eqref{eq1.15} we obtain the system of linear equations relative to $a(\lambda)$, $b(\lambda)$
$$\left\{
\begin{array}{lll}
{\displaystyle\alpha a(\lambda)W(\lambda)+\frac{2b(\lambda)}{\Delta(0,\lambda)}\left\{\alpha F(\lambda)+\alpha F^*(\lambda)+\alpha G(\lambda)-\frac{\Delta(0,\lambda)}2\right\}=0;}\\
\displaystyle{a(\lambda)\Delta(0,\lambda)+\frac{2b(\lambda)}{\Delta(0,\lambda)}N(\lambda)=0;}
\end{array}\right.$$
non-trivial solution to which exists if the determinant $\Delta(\alpha,\lambda)$ of this system vanishes, $\Delta(\alpha,\lambda)=0$, where
\begin{equation}
\Delta(\alpha,\lambda)\stackrel{\d}{=}\Delta(0,\lambda)+2\alpha\left\{\frac{W(\lambda)N(\lambda)}{\Delta(0,\lambda)}-F(\lambda)-F^*(\lambda)-G(\lambda)\right\}.\label{eq1.17}
\end{equation}
The function $\Delta(\alpha,\lambda)$ \eqref{eq1.17} is said to be the characteristic function of the operator $L$ \eqref{eq1.9} which coincides with $\Delta(0,\lambda)$ \eqref{eq1.3} when $\alpha=0$. Its zeroes $\lambda_k(\alpha)$ define the eigenvalues $z_k(\alpha)=\lambda_k^2(\alpha)$ of the operator $L$ \eqref{eq1.9}, i. e.,
\begin{equation}
\sigma(L)=\left\{z_k(\alpha)=\lambda_k^2(\alpha);\Delta\left(\alpha,\lambda_k(\alpha)\right)=0\right\}.\label{eq1.18}
\end{equation}
\vspace{5mm}

{\bf 1.3.} Let us simplify the formula for the character function $\Delta(\alpha,\lambda)$ \eqref{eq1.17}. Let
\begin{equation}
\widetilde{v}(\lambda)=\int\limits_0^\pi e^{-i\lambda x}v(x)dx\label{eq1.19}
\end{equation}
be the Fourier transform of the function $v(x)$. Then the function $W(\lambda)$ \eqref{eq1.14} equals
$$W(\lambda)=\frac1{2i\lambda}\left\{\int\limits_0^\pi e^{i\lambda x}\overline{v(x)}dx-\int\limits_0^\pi e^{-i\lambda x}\overline{v(x)}dx+e^{i\lambda\pi}\int\limits_0^\pi e^{-i\lambda x}\overline{v}(x)dx\right.$$
$$\left.-e^{-i\lambda\pi}\int\limits_0^\pi e^{i\lambda x}\overline{v}(x)dx\right\}$$
or, in terms of $\widetilde v(\lambda)$ \eqref{eq1.19},
$$W(\lambda)=\frac1{2i\lambda}\left\{\widetilde{v}^*(\lambda)\left(1-e^{-i\lambda\pi}\right)-\widetilde{v}^*(-\lambda)\left(1-e^{i\lambda\pi}\right)\right\}.$$
Similarly, $N(\lambda)$ \eqref{eq1.16} equals
$$N(\lambda)=\frac12\left\{\widetilde{v}(\lambda)\left(1-e^{i\lambda\pi}\right)+\widetilde{v}(-\lambda)\left(1-e^{-i\lambda\pi}\right)\right\}.$$
Using these expressions for $W(\lambda)$ and $N(\lambda)$, we obtain
$$\frac{W(\lambda)N(\lambda)}{\Delta(0,\lambda)}=\frac1{4i\lambda\Delta(0,\lambda)}\left\{\widetilde{v}(\lambda)\widetilde{v}^*(\lambda)\left(1-e^{i\lambda\pi}\right)\left(1-e^{-i\lambda\pi}\right)\right.$$
\begin{equation}
\begin{array}{ccc}
+\widetilde{v}(\lambda)\widetilde{v}^*(-\lambda)\left(1-e^{i\lambda\pi}\right)\left(1-e^{-i\lambda\pi}\right)e^{i\lambda\pi}-\widetilde{v}(-\lambda)\widetilde{v}^*(\lambda)\left(1-e^{i\lambda\pi}\right)\left(1-e^{-i\lambda\pi}\right)
\\
\displaystyle{\left.\times e^{-i\lambda\pi}-\widetilde{v}(-\lambda)\widetilde{v}^*(-\lambda)\left(1-e^{i\lambda\pi}\right)\left(1-e^{-i\lambda\pi}\right)\right\}=\frac1{4i\lambda}\left\{\widetilde{v}(\lambda)\widetilde{v}^*(\lambda)\right.}\label{eq1.20}
\end{array}
\end{equation}
$$\left.-\widetilde{v}(-\lambda)\widetilde{v}^*(-\lambda)+\widetilde{v}(\lambda)\widetilde{v}^*(-\lambda)e^{i\lambda\pi}-\widetilde{v}(-\lambda)\widetilde{v}^*(\lambda)e^{-i\lambda\pi}\right\}.$$
The function $G(\lambda)$ \eqref{eq1.14} equals
\begin{equation}
G(\lambda=\frac1{2i\lambda}\left(\Phi^*(\lambda)-\Phi(-\lambda)\right)\label{eq1.21}
\end{equation}
where
\begin{equation}
\begin{array}{ccc}
{\displaystyle\Phi(\lambda)\stackrel{\d}{=}\int\limits_0^\pi\int\limits_0^xe^{-i\lambda(x-t)}\overline{v}(t)dtv(x)dx=\int\limits_0^\pi\int\limits_0^xe^{-i\lambda\xi}\overline{v(x-\xi)}d\xi v(x)dx}\\
{\displaystyle=\int\limits_0^\pi e^{-i\lambda\xi}\int\limits_\xi^\pi\overline{v(x-\xi)}v(x)dx.}
\end{array}\label{eq1.22}
\end{equation}
Consequently, $\Phi(\lambda)$ is the Fourier transform,
\begin{equation}
\Phi(\lambda)=\int\limits_0^\pi e^{-i\lambda x}g(x)dx,\label{eq1.23}
\end{equation}
of the convolution
\begin{equation}
g(x)\stackrel{\d}{=}\int\limits_x^\pi\overline{v(t-x)}v(t)dt.\label{eq1.24}
\end{equation}

\begin{lemma}\label{l1.1}
The entire function $\Phi(\lambda)$ \eqref{eq1.22} of exponential type ($\leq\pi$) satisfies the relation
\begin{equation}
\Phi(\lambda)+\Phi^*(\lambda)=\widetilde{v}(\lambda)\widetilde v^*(\lambda)\label{eq1.25}
\end{equation}
where $\widetilde v(\lambda)$ is given by \eqref{eq1.19}. Any entire function with indicator diagram $[0,i\pi]$ that satisfies equality \eqref{eq1.25} coincides with $\Phi(\lambda)$ where $g(x)$ is given by \eqref{eq1.24}.
\end{lemma}

P r o o f. Equality \eqref{eq1.25} follows from the formula for integration by parts,
$$\Phi(\lambda)=\int\limits_0^\pi\int\limits_0^xe^{i\lambda t}\overline{v(t)}dtd_x\int\limits_0^xe^{-i\lambda s}v(s)ds=\widetilde v(\lambda)\widetilde v^*(\lambda)-$$
$$-\int\limits_0^\pi\int\limits_0^xe^{-i\lambda s}v(s)ds\cdot e^{i\lambda x}\overline{v(x)}dx=\widetilde{v}(\lambda)\widetilde{v}^*(\lambda)-\Phi^*(\lambda).$$
Let $F(\lambda)$ be an entire function of exponential type with indicator diagram $[0,2\pi]$, then the Paley -- Wiener theorem implies that
$$F(\lambda)=\int\limits_0^\pi e^{-i\lambda t}f(t)dt\quad(f\in L^2(0,\pi)).$$
If $F(\lambda)$ satisfies relation \eqref{eq1.25}, then the difference $G(\lambda)=F(\lambda)-\Phi(\lambda)$ ($\Phi(\lambda)$ is given by \eqref{eq1.23}) has the property $G(\lambda)=G^*(\lambda)$ which implies
$$0=\int\limits_0^\pi e^{-i\lambda t}(f(t)-g(t))dt-\int\limits_0^\pi e^{i\lambda t}\left(\overline{f(t)}-\overline{g(t)}\right)dt$$
$$=\int\limits_{-\pi}^\pi e^{-i\lambda t}\left\{(f(t)-g(t))\chi_+(t)-\left(\overline{f(-t)}-\overline{g(-t)}\right)\chi_-(t)\right\}dt$$
($\chi_\pm(t)$ are characteristic functions of the sets $\mathbb{R}_\pm$). This and  Parseval's identity yield that $f(t)=g(t)$. $\blacksquare$

\begin{lemma}\label{l1.2}
The identity
\begin{equation}
\begin{array}{ccc}
\displaystyle{F(\lambda)+F^*(\lambda)=\frac1{4i\lambda}\left\{(\Phi(\lambda)-\Phi(-\lambda))\left(e^{i\lambda\pi}+e^{-i\lambda\pi}\right)\right.}\\
+\left.\left(\widetilde{v}(\lambda)+\widetilde{v}(-\lambda)\right)\left(e^{i\lambda\pi}\widetilde{v}^*(-\lambda)-e^{-i\lambda\pi}\widetilde{v}^*(\lambda)\right)\right\}
\end{array}\label{eq1.26}
\end{equation}
is true for $F(\lambda)$ \eqref{eq1.14} where $\widetilde{v}(\lambda)$ and $\Phi(\lambda)$ are given by \eqref{eq1.19} and \eqref{eq1.22}.
\end{lemma}

P r o o f. Formula \eqref{eq1.14} implies
$$F(\lambda)=\frac1{4i\lambda}\int\limits_0^\pi\left(e^{i\lambda(\pi-x)}-e^{-i\lambda(\pi-x)}\right)\int\limits_0^x\left(e^{i\lambda t}+e^{-i\lambda t}\right)v(t)dt\overline{v}(x)dx$$
$$=\frac 1{4i\lambda}\left\{e^{i\lambda\pi}\int\limits_0^\pi\int\limits_0^xe^{i\lambda t}v(t)dte^{-i\lambda x}\overline{v(x)}dx+e^{i\lambda\pi}\int\limits_0^\pi\int\limits_0^xe^{-i\lambda t}v(t)dte^{-i\lambda x}\overline{v(x)}dx\right.$$
$$\left.-e^{-i\lambda\pi}\int\limits_0^\pi\int\limits_0^xe^{i\lambda t}v(t)dte^{i\lambda x}\overline{v}(x)dx-e^{-i\lambda\pi}\int\limits_0^\pi\int\limits_0^xe^{-i\lambda t}v(t)dte^{i\lambda x}\overline{v(x)}dx\right\}$$
$$=\frac1{4i\lambda}\left\{e^{i\lambda\pi}\Phi^*(-\lambda)-e^{-i\lambda\pi}\Phi^*(\lambda)+e^{i\lambda\pi}\psi(\lambda)-e^{-i\lambda\pi}\psi(-\lambda)\right\}$$
where
$$\psi(\lambda)=\int\limits_0^\pi\int\limits_0^xe^{-i\lambda t}v(t)dte^{-i\lambda x}\overline{v(x)}dx.$$
Similarly to \eqref{eq1.25}, the identity
$$\psi(\lambda)+\psi^*(\lambda)=\widetilde{v}(\lambda)\widetilde{v}^*(-\lambda)$$
is true for $\psi(\lambda)$. This identity and \eqref{eq1.25} implies
$$F(\lambda)+F^*(\lambda)=\frac1{4i\lambda}\left\{e^{i\lambda\pi}\Phi^*(-\lambda)-e^{-i\lambda\pi}\Phi^*(\lambda)+e^{i\lambda\pi}\psi(\lambda)-e^{-i\lambda\pi}\psi(-\lambda)\right.$$
$$\left.-e^{-i\lambda\pi}\Phi(-\lambda)+e^{i\lambda\pi}\Phi(\lambda)-e^{-i\lambda\pi}\psi^*(\lambda)+e^{i\lambda\pi}\psi^*(-\lambda)\right\}$$
$$=\frac1{4i\lambda}\left\{e^{i\lambda\pi}\Phi(\lambda)-e^{-i\lambda\pi}\Phi(-\lambda)+e^{i\lambda\pi}\left(\widetilde{v}(-\lambda)\widetilde{v}^*(-\lambda)-\Phi(-\lambda)\right)\right.$$
$$\left.-e^{-i\lambda\pi}\left(\widetilde{v}(\lambda)\widetilde{v}^*(\lambda)-\Phi(\lambda)\right)+e^{i\lambda\pi}\widetilde{v}(\lambda)\widetilde{v}^*(-\lambda)-e^{-i\lambda\pi}\widetilde{v}(-\lambda)\widetilde{v}^*(\lambda)
\right\}$$
$$=\frac1{4i\lambda}\left\{\Phi(\lambda)\left(e^{i\lambda\pi}+e^{-i\lambda\pi}\right)-\Phi(-\lambda)\left(e^{i\lambda\pi}+e^{-i\lambda\pi}\right)+\widetilde{v}(\lambda)\left(e^{i\lambda\pi}\widetilde{v}^*(-\lambda)\right.\right.$$
$$\left.\left.-e^{-i\lambda\pi}v^*(\lambda)\right)+\widetilde{v}(-\lambda)\left(e^{i\lambda\pi}\widetilde{v}^*(-\lambda)-e^{-i\lambda\pi}\widetilde{v}^*(\lambda)\right)\right\}$$
which proves \eqref{eq1.26}. $\blacksquare$

As a result, we arrive at the following statement.

\begin{theorem}\label{t1.1}
The representation
\begin{equation}
\Delta(\alpha,\lambda)=\Delta(0,\lambda)+\frac\alpha{2i\lambda}(R(\lambda)-R(-\lambda))\label{eq1.27}
\end{equation}
is true for the characteristic function \eqref{eq1.17}, where $R(\lambda)$ is expressed via $\widetilde{v}(\lambda)$ \eqref{eq1.19} and $\Phi(\lambda)$ \eqref{eq1.22} by the formula
\begin{equation}
R(\lambda)\stackrel{\d}{=}\left(1-e^{-i\lambda\pi}\right)\left\{\Phi(\lambda)\left(1-e^{i\lambda\pi}\right)-\widetilde{v}(\lambda)\widetilde{v}^*(\lambda)\right\}.\label{eq1.28}
\end{equation}
\end{theorem}

P r o o f. The relations \eqref{eq1.17}, \eqref{eq1.20}, \eqref{eq1.21}, \eqref{eq1.26} imply
$$\Delta(\alpha,\lambda)=\Delta(0,\lambda)+\frac\alpha{2i\lambda}\left\{\widetilde{v}(\lambda)\widetilde{v}^*(\lambda)-\widetilde{v}(-\lambda)\widetilde{v}^*(-\lambda)+e^{i\lambda\pi}\widetilde{v}(\lambda)\widetilde{v}^*(-\lambda)
\right.$$
$$-e^{-i\lambda\pi}\widetilde{v}(-\lambda)\widetilde{v}^*(\lambda)-(\Phi(\lambda)-\Phi(-\lambda))\left(e^{i\lambda\pi}+e^{-i\lambda\pi}\right)-\widetilde{v}(\lambda)e^{i\lambda\pi}\widetilde{v}^*(-\lambda)$$
$$+\widetilde{v}(\lambda)e^{-i\lambda\pi}\widetilde{v}^*(\lambda)-\widetilde{v}(-\lambda)e^{i\lambda\pi}\widetilde{v}^*(-\lambda)+\widetilde{v}(-\lambda)e^{-i\lambda\pi}\widetilde{v}^*(\lambda)$$
$$\left.-2\left[\widetilde{v}(\lambda)\widetilde{v}^*(\lambda)-\Phi(\lambda)-\widetilde{v}(-\lambda)\widetilde{v}^*(-\lambda)+\Phi(-\lambda)\right]\right\}=\Delta(0,\lambda)$$
$$+\frac\alpha{2i\lambda}\left\{\widetilde{v}(\lambda)\widetilde{v}^*(\lambda)\left(-1+e^{-i\lambda\pi}\right)-\widetilde{v}(-\lambda)\widetilde{v}^*(-\lambda)\left(-1+e^{i\lambda\pi}\right)\right.$$
$$\left.+(\Phi(\lambda)-\Phi(-\lambda))\left(1-e^{i\lambda\pi}\right)\left(1-e^{-i\lambda\pi}\right)\right\}.\blacksquare$$

Similarly to \eqref{eq1.4}, the function $\Delta(\alpha,\lambda)$ satisfies the relations
\begin{equation}
\Delta(\alpha,\lambda)=\Delta(\alpha,-\lambda);\quad\Delta^*(\alpha,\lambda)=\Delta(\alpha,\lambda).\label{eq1.29}
\end{equation}
The first equality is a corollary of formula \eqref{eq1.27}. To prove the second equality, we calculate
$$R^*(\lambda)=\left(1-e^{i\lambda\pi}\right)\left(1-e^{-i\lambda\pi}\right)\Phi^*(\lambda)-\left(1-e^{i\lambda\pi}\right)v^*(\lambda)v(\lambda)$$
$$=\left(1-e^{i\lambda\pi}\right)\left(1-e^{-i\lambda\pi}\right)\left(v^*(\lambda)v(\lambda)-\Phi(\lambda)\right)-\left(1-e^{i\lambda\pi}\right)v^*(\lambda)v(\lambda)$$
$$=-\left(1-e^{i\lambda\pi}\right)\left(1-e^{-i\lambda\pi}\right)\Phi(\lambda)+\left(1-e^{-i\lambda\pi}\right)v^*(\lambda)v(\lambda)$$
$$=-\left(1-e^{-i\lambda\pi}\right)\left\{\left(1-e^{i\lambda\pi}\right)\Phi(\lambda)-v^*(\lambda)v(\lambda)\right\}=-R(\lambda),$$
hence it follows that $\Delta^*(\alpha,\lambda)=\Delta(\alpha,\lambda)$.
\vspace{5mm}

{\bf 1.4.} Let us find the eigenfunctions of the operator $L$ \eqref{eq1.9}. We find $a(\lambda)$ from formula \eqref{eq1.15},
$$a(\lambda)=-\frac{2b(\lambda)}{\Delta(0,\lambda)}\frac{N(\lambda)}{\Delta(0,\lambda)},$$
then from \eqref{eq1.12} we obtain
\begin{equation}
\begin{array}{ccc}
{\displaystyle u(\lambda,x)=\frac{2b(\lambda)}{\Delta(0,\lambda)}\left\{-\frac{N(\lambda)}{\lambda\cdot\Delta(0,\lambda)}2\sin\frac{\lambda\pi}2\cos\lambda\left(\frac\pi2-x\right)+\int\limits_0^x\frac{\sin\lambda(x-t)}\lambda v(t)dt\right.}\\
{\displaystyle\left.+\frac{\sin\lambda(\pi-x)}\lambda\int\limits_0^x\cos\lambda tv(t)dt+\cos\lambda x\int\limits_x^\pi\frac{\sin\lambda(\pi-t)}\lambda v(t)dt\right\}.}
\end{array}\label{eq1.30}
\end{equation}
Since
$$N(\lambda)=\frac1i\sin\frac{\lambda\pi}2\left[\widetilde{v}(\lambda)e^{\frac{i\lambda\pi}2}-v(-\lambda)e^{-\frac{i\lambda\pi}2}\right]=2\sin\frac{\lambda\pi}2\int\limits_0^\pi\sin\lambda\left(\frac\pi 2-t\right)v(t)dt,$$
then
$$-\frac{N(\lambda)}{\lambda\cdot\Delta(0,\lambda)}2\sin\frac{\lambda\pi}2\cos\lambda\left(\frac\pi 2-x\right)=-\cos\lambda\left(\frac\pi 2-x\right)\int\limits_0^\pi\frac{\sin\lambda\left(\frac\pi 2-t\right)}\lambda v(t)dt,$$
and after elementary transformations from \eqref{eq1.30} we find the form of eigenfunctions of the operator $L$ \eqref{eq1.9}.

\begin{theorem}\label{t1.2}
To each zero $\lambda=\lambda_k(\alpha)\not\in\sigma\left(L_0\right)$ \eqref{eq1.6} of the characteristic function $\Delta(\alpha,\lambda)$ \eqref{eq1.27} there corresponds an eigenfunction of the operator $L$ \eqref{eq1.9},
\begin{equation}
u(\lambda,x)=\int\limits_0^x\cos\lambda\left(\frac\pi 2-x+t\right)v(t)dt+\int\limits_x^\pi\cos\lambda\left(\frac\pi 2-t+x\right)v(t)dt.\label{eq1.31}
\end{equation}
\end{theorem}

\section{Abstract formulation of the problem. Spectrum of the operator $L$}\label{s2}

{\bf 2.1.} Consider in a Hilbert space $H$ a self-adjoint operator $L_0$  with dense domain $\mathfrak{D}\left(L_0\right)$, spectrum of which consists of the eigenvalues $z_k$ of spectral multiplicity
\begin{equation}
\sigma\left(L_0\right)=\left\{z_k\in\mathbb{R};k\in\mathbb{Z}_+\right\}\label{eq2.1}
\end{equation}
assuming that numbers $z_k$ are numbered in increasing order and do not have finite limit points, i. e., $z_k\rightarrow\infty$ ($k\rightarrow\infty$). To an eigenvalue $z_k$ there corresponds an eigenspace $G_k$ of finite dimension, $\dim G_k=m_k<\infty$. Further, we assume that $m_k\leq m<\infty$ ($\forall k$). Spectral decomposition of the operator $L_0$ is
\begin{equation}
L_0=\sum\limits_{k=0}^\infty z_kE_k,\label{eq2.2}
\end{equation}
where $E_k$ are the orthoprojectors onto $G_k$ and $E_kE_s=\delta_{k,s}E_k$. Resolvent $R_{L_0}(z)=\left(L_0-zI\right)^{-1}$ of the operator $L_0$ \eqref{eq2.2} equals
\begin{equation}
R_{L_0}(z)=\sum\limits_k\frac{E_k}{z_k-z}.\label{eq2.3}
\end{equation}
Similarly to \eqref{eq1.9}, we define the self-adjoint operator
\begin{equation}
L=L_\alpha\stackrel{\d}{=}L_0+\alpha\langle\cdot,v\rangle v,\label{eq2.4}
\end{equation}
where $\alpha\in\mathbb{R}$, $v\in H$, and domains of the operators $L$ and $L_0$ coincide. Resolvent $R_L(z)=(L-zI)^{-1}$ of the operator $L$ is expressed via the resolvent $R_{L_0}(z)$ by the formula
\begin{equation}
R_L(z)f=R_{L_0}(z)f-\alpha\frac{\left\langle R_{L_0}(z)f,v\right\rangle}{1+\alpha\left\langle R_{L_0}(z)v,v\right\rangle}\cdot R_{L_0}(z)f\label{eq2.5}
\end{equation}
($f\in H$), and taking \eqref{eq2.3} into account we obtain
\begin{equation}
R_L(z)f=\sum\limits_k\frac{f_k}{z_k-z}-\frac{\displaystyle\alpha\sum\limits_s\frac{\left\langle f_s,v_s\right\rangle}{z_s-z}}{\displaystyle 1+\alpha\sum\limits_s\frac{\left\|v_s\right\|^2}{z_s-z}}\cdot\sum_k\frac{v_k}{z_k-z},\label{eq2.6}
\end{equation}
where $f_k=E_kf$, $v_k=E_kv$ are projections of $f$ and $v$ onto the subspace $G_k$. Residue of the function $R_L(z)f$ in the point $z=z_p$ equals
$$\res_{z_p}R_L(z)=f_p-\left\langle f,v_p\right\rangle\frac{v_p}{\left\|v_p\right\|^2}\stackrel{\d}{=}f_p^1.$$
Hence it follows: 1) if $v_p=0$, then the totality $f_p^1$ ($=f_p$) runs through all the subspace $G_p$, when $f\in H$, and consequently spectral multiplicity of $z_p$ equals $m_p$; 2) for $v_p\not=0$, the vectors $f_p^1$ belong to $G_p$ and are orthogonal to $v_p$, and thus the totality of $f_p^1$ forms the subspace $G_p^1$ in $G_p$ of dimension $m_p-1$. Therefore it is natural to distinguish two not-intersecting subsets in $\sigma\left(L_0\right)$ \eqref{eq2.1}:
\begin{equation}
\begin{array}{lll}
\sigma_0\stackrel{\d}{=}\left\{z_k\in\sigma\left(L_0\right):v_k=0\right\};\\
\sigma_1\stackrel{\d}{=}\left\{z_k\in\sigma\left(L_0\right):v_k\not=0,m_k>1\right\}.
\end{array}\label{eq2.7}
\end{equation}
It is obvious that $\sigma_0$ and $\sigma_1$ belong to the spectrum of the operator $L$ \eqref{eq2.4}, $\sigma_0$, $\sigma_1\subset\sigma(L)$, besides, spectral multiplicities of the points $z_k\in\sigma_0$ do not change and equal $m_k$, and for $z_k\in\sigma_1$ spectral multiplicities decrease by 1 and equal $m_k-1$. Formula \eqref{eq2.6} implies that zeroes of the function
\begin{equation}
Q(z)\stackrel{\d}{=}1+\alpha\left\langle R_{L_0}(z)v,v\right\rangle=1+\alpha\sum\limits_k\frac{\left\|v_k\right\|^2}{z_k-z}\label{eq2.8}
\end{equation}
belong to the spectrum $\sigma(L)$ of the operator $L$, besides, summation in \eqref{eq2.8} is carried out by such $k$ that $z_k\in\left(\sigma\left(L_0\right)\backslash\sigma_0\right)$. Function $Q(z)$ has simple zeroes $\mu_k\in\mathbb{R}$ interchanging with points from the set $\sigma\left(L_0\right)\backslash\sigma_0$. This set of zeroes $\left\{\mu_k\right\}$ we denote by $\sigma_2$,
\begin{equation}
\sigma_2\stackrel{\d}{=}\left\{\mu_k:Q\left(\mu_k\right)=0\right\}.\label{eq2.9}
\end{equation}
A number $\mu_k$ may coincide with some $z_s$ from $\sigma_0$, i. e., $\sigma_2\cap\sigma_0\not=\emptyset$, and in this case spectral multiplicity of $z_s$ increases by 1 and equals $m_s+1$.

\begin{theorem}\label{t2.1}
If spectrum $\sigma\left(L_0\right)$ \eqref{eq1.1} of a self-adjoint operator $L_0$ is discrete, besides, spectral multiplicities of the eigenvalues $z_k$ equal $m_k<m<\infty$ ($\forall k$) and $\sigma\left(L_0\right)$ does not have finite limit points, then spectrum $\sigma(L)$ of the operator $L$ \eqref{eq2.4} is the union of four non-intersecting sets,
\begin{equation}
\sigma(L)=\left\{\sigma_0\backslash\left(\sigma_0\cap\sigma_2\right)\right\}\cup\left\{\sigma_2\backslash\left(\sigma_0\cap\sigma_3\right)\right\}\cup\left\{\sigma_0\cap\sigma_2\right\}\cup\sigma_1,\label{eq2.10}
\end{equation}
besides, ${\rm (i)}$ spectral multiplicity $m_k$ of the points $z_k$ from $\sigma_0\backslash\left(\sigma_0\cap\sigma_2\right)$ remains unchanged  and $\sigma_0$ is given by \eqref{eq2.7};

${\rm(ii)}$ the points $\mu_k$ from $\sigma_2\backslash\left(\sigma_0\cap\sigma_2\right)$ have spectral multiplicity 1 and $\sigma_2$ is given by \eqref{eq2.9};

${\rm(iii)}$ the points $z_k$ from $\sigma_0\cap\sigma_2$, which coincide with zeroes $\mu_s$ of the function $Q(z)$ \eqref{eq2.8}, have spectral multiplicity $m_k+1$;

${\rm(iv)}$ the points $z_k\in\sigma_1$ have spectral multiplicity $m_k-1$;

${\rm(v)}$ points of the sets $\sigma_2$ and $\sigma\left(L_0\right)\backslash\sigma_0$ alternate.
\end{theorem}

\begin{remark}\label{r2.1}
If points of the specrum $z_k\in\sigma\left(L_0\right)$ \eqref{eq2.1} of the operator $L_0$ satisfy the separability property
\begin{equation}
\inf\limits_k\left(z_{k+1}-z_k\right)=d>0,\label{eq2.11}
\end{equation}
then the set $\sigma_0\cap\sigma_2$ is finite.
\end{remark}

If the set $\sigma_0\cap\sigma_2$ is infinite, then there exists the infinite sequence of roots $\mu_{k_p}$ of the function $Q(z)$ ($Q\left(\mu_{k_p}\right)=0$) coinciding with $z_p$ from $\sigma_0$. Formula \eqref{eq2.8} implies
$$1+\alpha\sum\limits_k\frac{\left\|v_k\right\|^2}{z_k-z_p}=0\quad(\forall z_p\in\sigma_0\cap\sigma_2).$$
Series in this sum converges uniformly by $p$ since $\left|z_k-z_p\right|>d$ and $\sum\limits_k\left\|v_k\right\|^2<\infty$ ($v\in H$). Proceeding to the limit as $p\rightarrow\infty$ and considering that $z_p\rightarrow\infty$, we obtain $0=1$, which is impossible.

\begin{remark}\label{r2.2}
Spectrum $\sigma\left(L_0\right)$ \eqref{eq1.6} of the operator $L_0$ \eqref{eq1.1}, \eqref{eq1.2} has separability property \eqref{eq2.11} and every eigenvalue $z_k(0)$ \eqref{eq1.6} has finite spectral multiplicity $m_k=2$ ($k\geq1$), $m_0=1$. Therefore Theorem \ref{t2.1} is true for $L$ \eqref{eq1.9}. For the convenience sake, we formulate it separately.
\end{remark}

\begin{theorem}\label{t2.2}
Spectrum $\sigma(L)$ of the operator $L$ \eqref{eq1.9} is the union of four non-intersecting sets,
\begin{equation}
\sigma(L)=\left\{\sigma_0\backslash\left(\sigma_0\cap\sigma_2\right)\right\}\cup\left\{\sigma_2\backslash\left(\sigma_0\cap\sigma_2\right)\right\}\cup\left\{\sigma_0\cap\sigma_2\right\}\cup\sigma_1.\label{eq2.12}
\end{equation}
Besides,

${\rm(i)}$ spectral multiplicity ($m_0=0$, $m_k=2$, $k\geq1$) of the points $z_k(\alpha)=z_k(0)$ from $\sigma_0\backslash\left(\sigma_0\cap\sigma_2\right)$ remains unchanged;

$\rm{(ii)}$ the points $z_k(\alpha)=\mu_k\in\sigma_2\backslash(\sigma_0\cap\sigma_2)$ have multiplicity $1$;

$\rm{(iii)}$ the set $\sigma_0\cap\sigma_2$ is finite and its points $z_k(\alpha)=z_k(0)=\mu_p\in\sigma_0\cap\sigma_2$ have multiplicity $m_k=3$ ($k\geq1$) and $m_0=2$ ($k=0$);

$\rm{(iv)}$ the points $z_k(\alpha)=z_k(0)\in\sigma_1$ \eqref{eq2.7} have multiplicity $1$;

$\rm{(v)}$ the points $\mu_k$ from $\sigma_2$ \eqref{eq2.9} alternate with the points $z_k(0)$ from $\sigma\left(L_0\right)\backslash\sigma_0$.
\end{theorem}

This theorem implies the description of eigenfunctions of the operator $L$ \eqref{eq1.9}.

\begin{theorem}\label{t2.3}
In accordance with decomposition \eqref{eq2.12} of the spectrum $\sigma(L)$, eigenfunctions of the operator $L$ \eqref{eq1.9} are as follows:

$\rm{(i)}$ to the points $z_k(\alpha)=z_k(0)\in\sigma_0\backslash\left(\sigma_0\cap\sigma_2\right)$ there correspond the two-dimensional eigenspaces $G_k$ \eqref{eq1.7}, for $z_k(0)=4k^2$ ($k\geq1$), and the one-dimensional eigenspace $\displaystyle{\left\{\xi\frac1{\sqrt{\pi}}\right\}}$, for $z_k(0)=0$;

$\rm{(ii)}$ to the eigenvalues $z_k(\alpha)=\mu_k\in\sigma_2\backslash\left(\sigma_0\cap\sigma_2\right)$ there correspond the eigenfunctions $u\left(\sqrt{\mu_k},x\right)$ \eqref{eq1.31}, where $\mu_k$ are zeroes \eqref{eq2.9} of the function $Q(z)$ \eqref{eq2.8};

$\rm{(iii)}$ to the points $z_k(\alpha)=z_k(0)=\mu_p\in\sigma_0\cap\sigma_2$ there corresponds the eigenspace $\span\left\{E_k+u\left(\sqrt{\mu_p},x\right)\right\}$;

$\rm{(iv)}$ to the eigenvalues $z_k(0)\in\sigma_1$ there correspond the eigenfunctions ${\displaystyle\left(c_k^2+s_k^2\right)^{-1/2}}$ $\displaystyle{\left(s_k\sqrt{\frac2\pi}\cos 2kx-s_k\sqrt{\frac2\pi}\sin 2kx\right)}$, where $s_k=\displaystyle{\left\langle v,\sqrt{\frac2\pi}\sin 2kx\right\rangle}$,\\ $\displaystyle{c_k=\left\langle v,\sqrt{\frac2\pi}\cos 2kx\right\rangle}$.
\end{theorem}
\vspace{5mm}

{\bf 2.2.} Let us study in detail the function $Q(z)$ \eqref{eq2.8} characterizing change of the spectrum of the operator $L_0$ under the one-dimensional perturbation \eqref{eq2.4}. We express the function $Q\left(\lambda^2\right)$ via the characteristic functions $\Delta(\alpha,\lambda)$ and $\Delta(0,\lambda)$.

If $R_{L_0}\left(\lambda^2\right)f=y$, then $L_0y-\lambda^2y=f$, this coincides with equation \eqref{eq1.10}, the solution to which is given by \eqref{eq1.11}. We find $a(\lambda)$ from the second boundary condition for $y(\lambda,x)$ and after some elementary transformations we obtain that
\begin{equation}
\begin{array}{ccc}
{\displaystyle y(\lambda,x)=\left(R_{L_0}\left(\lambda^2\right)f\right)(x)=-\frac1{\displaystyle2\lambda\sin\frac{\pi\lambda}2}\left\{\int\limits_0^x\cos\lambda\left(\frac\pi2-x+t\right)f(t)dt\right.}\\
\displaystyle{\left.+\int\limits_x^\pi\cos\lambda\left(\frac\pi2-t+x\right)f(t)dt\right\}.}
\end{array}\label{eq2.13}
\end{equation}
Hence
$$\left\langle R_{L_0}\left(\lambda^2\right)v,v\right\rangle=-\frac1{\displaystyle 2\lambda\sin\frac{\pi\lambda}2}\left\{\int\limits_0^\pi\int\limits_0^x\cos\lambda\left(\frac\pi2-x+t\right)v(t)dt\overline{v(x)}dx\right.$$
$$\left.+\int\limits_0^\pi\int\limits_x^\pi\cos\lambda\left(\frac\pi2-t+x\right)v(t)dt\overline{v(x)}dx\right\}=-\frac1{\displaystyle 4\lambda\sin\frac{\pi\lambda}2}\left\{e^{i\lambda\frac\pi2}\int\limits_0^\pi\int\limits_0^xe^{i\lambda t}v(t)dt\right.$$
$$\times e^{-i\lambda x}\overline{v(x)}dx+e^{-i\lambda\frac\pi2}\int\limits_0^\pi\int\limits_0^xe^{-i\lambda t}v(t)dte^{i\lambda x}\overline{v(x)}dx+e^{i\lambda\frac\pi2}\int\limits_0^\pi\int\limits_x^\pi e^{-i\lambda t}v(t)dt$$
$$\left.\times e^{i\lambda x}\overline{v(x)}dx+e^{-i\lambda\frac\pi2}\int\limits_0^\pi\int\limits_x^\pi e^{i\lambda t}v(t)dte^{-i\lambda x}\overline{v(x)}dx\right\}=$$
$$=-\frac1{\displaystyle 4\lambda\sin\frac{\pi\lambda}2}\left\{e^{i\lambda\frac\pi2}\Phi^*(-\lambda)+e^{-\frac{i\lambda\pi}2}\Phi^*(\lambda)+e^{i\lambda\frac\pi2}\Phi(\lambda)+e^{-\frac{i\lambda\pi}2}\Phi(-\lambda)\right\}$$
in view of definition \eqref{eq1.21} of the function $\Phi(\lambda)$. Relation \eqref{eq1.25} yields
$$\left\langle R_{L_0}\left(\lambda^2\right)v,v\right\rangle=-\frac1{\displaystyle4\lambda\sin\frac{\lambda\pi}2}\left\{e^{i\lambda\frac\pi2}\left(\widetilde{v}(-\lambda)\widetilde{v}^*(-\lambda)-\Phi(-\lambda)\right)\right.$$
$$\left.+e^{-\frac{i\lambda\pi}2}\left(\widetilde{v}(\lambda)\widetilde{v}^*(\lambda)-\Phi(\lambda)\right)+e^{i\lambda\frac\pi2}\Phi(\lambda)+e^{-\frac{i\lambda\pi}2}\Phi(-\lambda)\right\}$$
$$=-\frac1{\displaystyle 4\lambda\sin\frac{\lambda\pi}2}\left\{e^{\frac{i\lambda\pi}2}\left[\widetilde{v}(-\lambda)\widetilde{v}^*(-\lambda)-\left(1-e^{-i\lambda\pi}\right)\Phi(-\lambda)\right]\right.$$
$$\left.+e^{-\frac{i\lambda\pi}2}\left[\widetilde{v}(\lambda)\widetilde{v}^*(\lambda)-\left(1-e^{i\lambda\pi}\right)\Phi(\lambda)\right]\right\}.$$
Taking \eqref{eq1.28} into account, we obtain
$$1+\alpha\left\langle R_{L_0}\left(\lambda^2\right)v,v\right\rangle=1-\frac\alpha{\displaystyle4\lambda\sin\frac{\lambda\pi}2}\left\{-\frac{e^{-\frac{i\lambda\pi}2}}{1-e^{-i\lambda\pi}}R(\lambda)-\frac{e^{\frac{i\lambda\pi}2}}{1-e^{i\lambda\pi}}
R(-\lambda)\right\}$$
$$=1+\frac\alpha{\displaystyle 4\lambda\sin\frac{\lambda\pi}2}\left\{\frac1{\displaystyle 2i\sin\frac{\lambda\pi}2}R(\lambda)-\frac1{\displaystyle 2i\sin\frac{\lambda\pi}2}R(-\lambda)\right\}$$
$$=1+\frac\alpha{2i\lambda\Delta(0,\lambda)}\{R(\lambda)-R(-\lambda)\}=\frac{\Delta(\alpha,\lambda)}{\Delta(0,\lambda)},$$
in view of \eqref{eq1.27}.

\begin{theorem}\label{t2.3}
The function $Q\left(\lambda^2\right)$ \eqref{eq2.8} is expressed via the characteristic function $\Delta(\alpha,\lambda)$ \eqref{eq1.27} of the operator $L$ \eqref{eq1.9} by the formula
\begin{equation}
Q\left(\lambda^2\right)=\frac{\Delta(\alpha,\lambda)}{\Delta(0,\lambda)}\label{eq2.14}
\end{equation}
where $\Delta(0,\lambda)$ is the characteristic function \eqref{eq1.3} of the operator $L_0$ \eqref{eq1.1}, \eqref{eq1.2}.
\end{theorem}

\section{Structure of the function $Q(z)$}\label{s3}

{\bf 3.1.} In order to obtain the multiplicative expansion of the characteristic function $\Delta(\alpha,\lambda)$ \eqref{eq1.27}, we use the well-known statements. Let us recall \cite{9} that a function $f(\lambda)$ is of class $C$ (Cartwright) if

(a) $f(\lambda)$ is an entire function of exponential type;

(b) the integral
\begin{equation}
\int\limits_{\mathbb{R}}\frac{\ln^+|f(x)|}{1+x^2}dx<\infty\label{eq3.1}
\end{equation}
converges. The function $\Delta(\alpha,\lambda)$ \eqref{eq1.27} is of class $C$. So the Paley -- Wiener theorem for $\widetilde{v}(\lambda)$ \eqref{eq1.19} and $\Phi(\lambda)$ \eqref{eq1.22} implies that $\Delta(\alpha,\lambda)$ satisfies condition (a). Convergence of integral \eqref{eq3.1} follows from analyticity and boundedness of the function $R(\lambda)$ \eqref{eq1.28}.

Indicator diagram \cite{9} of the function $\Delta(\alpha,\lambda)$ coincides with the segment of imaginary axis $\displaystyle{\left[-ih_\Delta\left(-\frac\pi2\right),ih_\Delta\left(\frac\pi2\right)\right]}$, besides, $\displaystyle{h_\Delta\left(-\frac\pi2\right)=h_\Delta\left(\frac\pi2\right)}$ in view of $|\Delta(\alpha,iy)|=|\Delta(\alpha,-iy)|$ \eqref{eq1.29}. The next statement follows from the more general Cartwright -- Levinson theorem \cite{9}.

\begin{theorem}\label{t3.1}
\cite{9}. For any function $f(\lambda)$ of class $C$ such that $f(0)\not=0$ and ${\displaystyle h_\Delta\left(-\frac\pi2\right)=h_\Delta\left(\frac\pi2\right)}$,
$$f(\lambda)=f(0)\lim\limits_{R\rightarrow\infty}\prod\limits_{\left|a_k\right|<R}\left(1-\frac\lambda{a_k}\right),$$
where $\left\{a_k\right\}_1^\infty$ are roots of the function $f(\lambda)$.
\end{theorem}

Relations \eqref{eq1.27}, \eqref{eq1.28} imply that $\Delta(\alpha,0)=-\alpha\pi\left|\widetilde v(0)\right|^2\not=0$. Applying Theorem \ref{t3.1} and taking \eqref{eq1.29} into account, we obtain
\begin{equation}
\Delta(\alpha,\lambda)=A\prod\limits_k\left(1-\frac{\lambda^2}{z_k(\alpha)}\right)\label{eq3.2}
\end{equation}
where $z_k(\alpha)=\lambda_k^2(\alpha)$ and $\lambda_k(\alpha)$ are zeroes of $\Delta(\alpha,\lambda)$. Equalities \eqref{eq2.8}, \eqref{eq2.14} yield
$$1+\alpha\sum\limits_k\frac{\left\|v_k\right\|^2}{z_k(0)-z}=\frac{\displaystyle A\prod\limits_k\left(1-\frac z{z_k(\alpha)}\right)}{\displaystyle 4\sin^2\sqrt{z}\frac\pi2},$$
and since \cite{9}
$$\sin\xi=\xi\prod\limits_k\left(1-\frac{\xi^2}{\pi^2k^2}\right);\quad 4\sin^2\sqrt{z}\frac\pi2=\pi^2z\prod\limits_k\left(1-\frac z{z_k(0)}\right),$$
then
$$1+\alpha\sum\limits_k\frac{\left\|v_k\right\|^2}{z_k(0)-z}=\frac{\displaystyle A\prod\limits_k\left(1-\frac z{z_k(\alpha)}\right)}{\displaystyle\pi^2z\prod\limits_k\left(1-\frac z{z_k(0)}\right)^2}.$$
Using Theorem \ref{t2.2}, we obtain
\begin{equation}
1+\alpha\sum\limits_k\frac{\left\|v_k\right\|^2}{z_k(0)-z}=\frac{\displaystyle A\prod\limits_k\left(1-\frac z{\mu_k}\right)}{\displaystyle\pi^2z\prod\limits_k\left(1-\frac z{z_k(0)}\right)}\label{eq3.3}
\end{equation}
where summation in the left-hand side of this equality and product in the denominator in the right-hand side is carried out only by $z_k(0)$ belonging to $\sigma\left(L_0\right)\backslash\sigma_0$. We obtain the number $A$ from \eqref{eq3.3},
\begin{equation}
\frac1A=\lim\limits_{y\rightarrow\infty}\frac{\displaystyle\prod\limits_k\left(1-\frac{iy}{\mu_k}\right)}{\displaystyle\pi^2iy\prod\limits_k\left(1-\frac{iy}{z_k(0)}\right)}.\label{eq3.4}
\end{equation}
Calculating the residue at the point $z_p$ in equality \eqref{eq3.3}, we obtain
\begin{equation}
\begin{array}{lll}
{\displaystyle\alpha\left\|v_p\right\|^2=A\frac{\mu_p-z_p(0)}{\pi^2\mu_p}\prod\limits_{k\not=p}\frac{z_k(0)}{\mu_k}\left(1-\frac{z_k(0)-\mu_k}{z_k(0)-z_p(0)}\right)\quad(p\not=0);}\\
{\displaystyle\alpha\left\|v_0\right\|^2=-\frac A{\pi^2}\quad(p=0).}
\end{array}\label{eq3.5}\end{equation}

\begin{theorem}\label{t3.2}
The numbers $\alpha\left\|v_p\right\|^2$ ($p\in\mathbb{Z}_+$), where $v_p=E_pv$ is the projection of the function $v$ onto the subspace $G_p$ \eqref{eq1.7}, are obtained unambiguously from the spectrum $\sigma(L)$ \eqref{eq2.12} of the operator $L$ \eqref{eq1.9} using the formulas \eqref{eq3.5}.
\end{theorem}

\begin{remark}\label{r3.1}
The number $\alpha$ is calculated from the equalities \eqref{eq3.5} unambiguously if we initially consider $\|v\|=1$, i. e., $\sum\left\|v_p\right\|^2=1$ (cf. Remark \ref{r1.1}). Consequently, the numbers $\alpha$ and $\left\|v_p\right\|^2$ are defined from $\sigma(L)$ unambiguously.

In the abstract formulation of the problem, when $L_0$ is not necessarily defined by \eqref{eq1.1}, to obtain $\alpha$ and $\left\|v_p\right\|^2$ we must know the two spectra $\sigma\left(L_0\right)$ and $\sigma(L)$.
\end{remark}
\vspace{5mm}

{\bf 3.2.} Let us present another method of calculation of numbers $\alpha\left\|v_p\right\|^2$. To do this, we write equality \eqref{eq2.14} as
\begin{equation}
1+\alpha\sum\limits_k\frac{\left\|v_k\right\|^2}{z_k(0)-z}=\frac{\Delta(\alpha,\lambda)}{\Delta(0,\lambda)}\quad\left(z=\lambda^2\right),\label{eq3.6}
\end{equation}
where summation in \eqref{eq3.6} is carried out by those $k$ for which $z_k(0)\in\sigma\left(L_0\right)\backslash\sigma_0$, then
\begin{equation}
\alpha\left\|v_p\right\|^2=\lim\limits_{\lambda\rightarrow\lambda_p(0)}\frac{\left(\lambda_p^2(0)-\lambda^2\right)\Delta(\alpha,\lambda)}{\Delta(0,\lambda)}.\label{eq3.7}
\end{equation}
Since $\lambda_p(0)\not=0$ is the second order zero of $\Delta(0,\lambda)$ and the first order zero of $\Delta(\alpha,\lambda)$, then taking into account that
$$\frac{\left(\lambda_p^2(0)-\lambda^2\right)\Delta(\alpha,\lambda)}{\Delta(0,\lambda)}=-\left(\lambda_p(0)+\lambda\right)\frac{\displaystyle\frac{\Delta(\alpha,\lambda)-\Delta\left(\alpha,\lambda_p(0)\right)}{\lambda-\lambda_p(0)}}
{\displaystyle\frac{\Delta(0,\lambda)-\Delta\left(\alpha,\lambda_p(0)\right)}{\left(\lambda-\lambda_p(0)\right)^2}},$$
from \eqref{eq3.7} we obtain
$$\alpha\left\|v_p\right\|^2=-4\lambda_p(0)\frac{\Delta'\left(0,\lambda_p(0)\right)}{\Delta''\left(0,\lambda_p(0)\right)}.$$
Hence from $\Delta''(0,\lambda)=2\pi^2\cos\lambda\pi$ we find
$$\alpha\left\|v_p\right\|^2=-\frac2{\pi^2}\lambda_p(0)\frac{\Delta'\left(0,\lambda_p(0)\right)}{\cos\lambda_p(0)\pi}\quad(\lambda_p=0)$$
or
\begin{equation}
\alpha\left\|v_p\right\|^2=-\frac{4p}{\pi^2}\Delta'(\alpha,2p)\quad(p\not=0).\label{eq3.8}
\end{equation}
Similarly,
\begin{equation}
\alpha\left|v_0\right|^2=-\frac1{\pi^2}\Delta(\alpha,0)\quad(p=0).\label{eq3.9}
\end{equation}

\begin{theorem}\label{t3.3}
The numbers $\alpha\left\|v_p\right\|^2$ are calculated via the characteristic function $\Delta(\alpha,\lambda)$ \eqref{eq1.27} by the formulas \eqref{eq3.8}, \eqref{eq3.9}.
\end{theorem}

\section{Inverse problem. Description of inverse problem data}\label{s4}

{\bf 4.1.} Theorem \ref{t3.2} and Remark \ref{r3.1} imply that the number $\alpha$ and squares of the norms $\left\|v_k\right\|^2$ of the projections $v_k=E_kv$ of the function $v$ onto the spaces $G_k$ \eqref{eq1.7} are unambiguously found from the spectrum $\sigma(L)$ of the operator $L$ \eqref{eq1.9}. Due to arbitrariness of an orthonormal basis in $G_k$, this means that $v(x)$ is recovered from $\left\|v_k\right\|^2$ with a large degree of uncertainty.

Since
\begin{equation}
v_k=c_k(v)\sqrt{\frac2\pi}\cos2kx+s_k(v)\sqrt\frac2\pi\sin2kx\label{eq4.1}
\end{equation}
where
\begin{equation}
c_k(v)=\left\langle v(x),\sqrt\frac2\pi\cos 2kx\right\rangle,\quad s_k(v)=\left\langle v(x),\sqrt2\pi\sin 2kx\right\rangle\label{eq4.2}
\end{equation}
are Fourier coefficients of the function $v(x)$, then
\begin{equation}
\left\|v_k\right\|^2=\left|c_k(v)\right|^2+\left|s_k(v)\right|^2.\label{eq4.3}
\end{equation}
We represent $L^2(0,\pi)$ as an orthogonal sum of subspaces,
\begin{equation}
L^2(0,\pi)=L_+\oplus L_-\label{eq4.4}
\end{equation}
where
\begin{equation}
L_{\pm}=\left\{f_\pm(x)=\frac12(f(x)\pm f(\pi-x)):f\in L^2(0,\pi)\right\}.\label{eq4.5}
\end{equation}
Orthogonality $L_+\perp L_-$ is a corollary of equalities
$$f_+(\pi-x)=f_+(x),\quad f_-(\pi-x)=-f_-(x).$$

\begin{remark}\label{r4.1} Since $\cos 2kx\in L_+$ and $\sin 2kx\in L_-$, the Fourier coefficients \eqref{eq4.2} of the function $f_\pm$ are
\begin{equation}
\begin{array}{lll}
c_k(f_+)=c_k(f),\quad s_k(f_+)=0;\\
c_k(f_-)=0,\quad s_k(f_-)=s_k(f).
\end{array}\label{eq4.6}
\end{equation}
\end{remark}

Theorem \ref{t3.2} and Remark \ref{r4.1} implies the following statement.

\begin{theorem}\label{t4.1}
The numbers $\alpha c_k(v)|^2$ and $\alpha|s_k(v)|^2$ are unambiguously found from the two spectra $\sigma\left(L_\pm\right)$ of operators $L_\pm=L\left(\alpha,v_\pm\right)$ \eqref{eq1.9} where $v_\pm=P_\pm v$, $P_\pm$ are the orthoprojectors onto the subspaces $L_\pm$ \eqref{eq4.5}.
\end{theorem}

Number $\alpha$ and coefficients $|c_k(v)|$ and $|s_k(v)|$ are calculated from $\alpha|c_k(v)|^2$ and $\alpha|s_k(v)|^2$ under the condition $\|v\|=1$ (cf. Remark \ref{r3.1}). This allows one to (ambiguously) recover $v$. Formula \eqref{eq4.3} implies that instead of spectra $\sigma\left(L_\pm\right)$ we can take spectra $\sigma(L)$ and $\sigma\left(L_+\right)$ (or $\sigma(L)$ and $\sigma\left(L_-\right)$).
\vspace{5mm}

{\bf 4.2.} In the case of real $v$, one can unambiguously recover $v(x)$. Formula \eqref{eq4.2} for a real function $v$ implies
\begin{equation}
\left\|v_k\right\|^2=c_k^2(v)+s_k^2(v).\label{eq4.7}
\end{equation}
Consider the function
\begin{equation}
w(x)=v(x)+x-\frac\pi2,\label{eq4.8}
\end{equation}
and since ${\displaystyle x-\frac\pi2\in L_-}$ \eqref{eq4.5}, then
$$c_k\left(x-\frac\pi2\right)=0,\quad s_k\left(x-\frac\pi2\right)=-\sqrt\frac\pi2\frac1k,$$
therefore
$$c_k(w)=c_k(v),\quad s_k(w)=s_k(v)-\sqrt\frac\pi2\frac1k,$$
and thus
$$\left\|w_k\right\|^2=c_k^2(v)+s_k^2(v)-\frac{\sqrt{2\pi}}ks_k(v)+\frac\pi{2k^2}.$$
Hence from \eqref{eq4.7} we find
\begin{equation}
s_k(v)=\frac k{\sqrt{2\pi}}\left\{\left\|v_k\right\|^2-\left\|w_k\right\|^2+\frac\pi{2k^2}\right\}.\label{eq4.9}
\end{equation}
Similarly, consider
\begin{equation}
\widehat{w}(x)=v(x)+\left(x-\frac\pi2\right)^2,\label{eq4.10}
\end{equation}
and taking ${\displaystyle\left(x-\frac\pi2\right)^2\in L_+}$ \eqref{eq4.5} we obtain
$$s_k\left(\left(x-\frac\pi2\right)^2\right)=0;\quad c_k\left(\left(x-\frac\pi2\right)^2\right)=\sqrt\frac\pi2\frac1{k^2}\quad(k\geq1);$$
$$c_0\left(\left(x-\frac\pi2\right)^2\right)=\frac{\pi^{5/2}}{12},$$
and thus
$$s_k(\widehat w)=s_k(v);\quad c_k(\widehat w)=c_k(v)+\sqrt{\frac\pi2}\frac1{k^2}\quad(k\geq1);$$
$$c_0(\widehat w)=c_0(v)+\frac{\pi^{5/2}}{12};$$
therefore
$$\left\|\widehat w_k\right\|^2=c_k^2(v)+\frac{\sqrt{2\pi}}{k^2}c_k(v)+\frac\pi{2k^4}+s_k^2(v)\quad(k\geq1);$$
$$\left\|\widehat w_0\right\|^2=c_0^2(v)+\frac{\pi^{5/2}}6c_0(v)+\frac{\pi^5}{144}\quad(k=0),$$
and according to \eqref{eq4.7}
\begin{equation}
\begin{array}{lll}
{\displaystyle c_k(v)=\frac{k^2}{\sqrt{2\pi}}\left\{\left\|\widehat w_k\right\|^2-\left\|v_k\right\|^2-\frac\pi{2k^4}\right\}\quad(k\geq1);}\\
{\displaystyle c_0(v)=\frac6{\pi^{5/2}}\left\{\left\|\widehat w_0\right\|^2-\left\|v_0\right\|^2-\frac{\pi^5}{144}\right\}\quad(k=0).}
\end{array}\label{eq4.11}
\end{equation}
Using Theorem \ref{t4.2}, we arrive at the statement.

\begin{theorem}\label{t4.2}
If $v(x)$ is real, then from the three spectra:

{\rm(i)} spectrum $\sigma(L)$ \eqref{eq2.12} of the operator $L=L(\alpha,v)$ \eqref{eq1.9};

{\rm(ii)} spectrum $\sigma\left(L_w\right)$ of the operator $L_w=L(\alpha,w)$ \eqref{eq1.9} where $w$ is given by \eqref{eq4.8};

{\rm(iii)} spectrum $\sigma\left(L_{\widehat w}\right)$ of the operator $L_{\widehat w}=L(\alpha,\widehat w)$ \eqref{eq1.9} where $\widehat w$ is given by \eqref{eq4.10};\\
the numbers $\alpha s_k(v)$ and $\alpha c_k(v)$ are unambiguously found via the formulas \eqref{eq3.5} and \eqref{eq4.9}, \eqref{eq4.11}.
\end{theorem}

\begin{remark}\label{r4.2}
Due to the normalization condition $\|v\|=1$, $\alpha^2$ is determined by the numbers $\alpha s_k(v)$, $\alpha c_k(v)$, and thus the number $\alpha$ is determined up to sign which is found from the pattern of intermittency of the numbers $\mu_k$ and $z_k(0)\in\sigma(L_0)\backslash\sigma_0$ (Theorem \ref{t2.2}, (v)). Thus, the number $\alpha$ and the real function $v(x)$ are unambiguously recovered from the three spectra $\sigma(L)$, $\sigma(L_w)$, $\sigma(L_{\widehat w})$.
\end{remark}

Instead of the functions ${\displaystyle\left(x-\frac\pi2\right)\in L_-}$ and ${\displaystyle\left(x-\frac\pi2\right)^2\in L_+}$, one can take any other functions $f$ and $g$ from the subspaces $L_-$ and $L_+$ \eqref{eq4.5} with non-zero Fourier coefficients $s_k(f)$ and $c_k(g)$ correspondingly.
\vspace{5mm}

{\bf 3.4.} Let us proceed to the definition of spectral data of the operator $L$ \eqref{eq1.9}. First, we give the formal definition of the set $\sigma(L)$ \eqref{eq2.12}.

\begin{definition}\label{d4.1}
Let $\sigma\left(L_0\right)$ be a countable set \eqref{eq1.6} with each point of multiplicity $2$. We chose a discrete subset $A$ from $\sigma(L_0)$ and no more then countable set of points $B$ which intermit with the elements of $\sigma(L_0)\backslash A.$ We match every such pair $A$, $B$ with the countable set
\begin{equation}
\sigma\left(L_0,A,B\right)\stackrel{\rm def}{=}\{A\backslash(A\cap B)\}\cup\{B\backslash(A\cap B)\}\cup\{A\cap B\}+\sigma_1\label{eq4.12}
\end{equation}
where $\sigma_1=\sigma(L_0)\backslash\{A\cup\{z(0)\}\}.$ Becides,

{\rm(i)} points $z_k(0)\in A\backslash\{A\cap B\}$ are of multiplicity $2$;

{\rm(ii)} points $\mu_k\in B\backslash\{A\cap B\}$ are of multiplicity $1$;

{\rm(iii)} the set $A\cap B$ is countable and its elements are of multiplicity $3$;

{\rm(iv)} points $z_k(0)\in\sigma_1$ are of multiplicity $1$.
\end{definition}

It is more convenient to define the numerical sequences $z_k(\alpha)$ which make up the spectrum $\sigma(L)$ of the operator $L$ \eqref{eq1.9} in terms of characteristic functions $\Delta(\alpha,\lambda)$ \eqref{eq1.26}. Let us proceed to definition of the class of functions which are characteristic for the operators of the type $L$ \eqref{eq1.9}.

\begin{definition}\label{d4.2}
Function $\Delta(\lambda)$ is of {\bf class} $JP\{z_k(0)\}$, where $z_k(0)=4k^2$ ($k\in\mathbb{Z}_+$) if

{\rm (i)} $\Delta(\lambda)$ is an entire function of exponential type $\sigma$ ($\leq2\pi$) such that $|\Delta(\lambda)|<\infty$ as $\lambda\in\mathbb{R}$ and $\Delta(-\lambda)=\Delta(\lambda)$, $\Delta^*(\lambda)=\Delta(\lambda)$;

{\rm (ii)} zeroes $\left\{\pm a_k\right\}$ of the function $\Delta(\lambda)$ ($a_k\not=0$ $\forall k$) have the property that $a_k^2\in\sigma(L_0,A,B)$ \eqref{eq4.12} for some $A$ and $B$;

{\rm (iii)} the function
\begin{equation}
F(z)=\frac{\Delta(\lambda)}{\Delta(0,\lambda)}\quad(\lambda^2=z)\label{eq4.13}
\end{equation}
($\Delta(0,\lambda)=2(1-\cos\lambda\pi)$ \eqref{eq1.3}) has the property
\begin{equation}
\lim\limits_{y\rightarrow\infty}F(iy)=1;\label{eq4.14}
\end{equation}
{\rm (iv)} the limit
\begin{equation}
\lim\limits_{y\rightarrow\infty}y|F(iy)-1|<\infty\label{eq4.15}
\end{equation}
exists.
\end{definition}

The function $\Delta(\alpha,\lambda)$ \eqref{eq1.26} is of class $JP\{z_k(0)\}$ due to \eqref{eq1.8} and Theorem \ref{t2.2}. Condition \eqref{eq4.15} in view of \eqref{eq3.6} provides convergence of the series ${\displaystyle \alpha\sum\|v_k\|^2<\infty}$, i. e., the fact that $v(x)$ belongs to the space $L^2(0,\pi)$.

\begin{theorem}\label{t4.3}
If $\Delta(\lambda)\in JP\{z_k(0)\}$, then there exists the operator $L$ \eqref{eq1.9}, characteristic function $\Delta(\alpha,\lambda)$ \eqref{eq1.26} of which coincides with $\Delta(\lambda)$.
\end{theorem}

P r o o f. Boundedness of $\Delta(\lambda)$ as $\lambda\in\mathbb{R}$ implies that the function $\Delta(\lambda)$ is of Cartwright class $C$. Using Theorem \ref{t3.1} and structure of the set $\sigma(L_0,A,B)$, we obtain that
\begin{equation}
F(z)=\frac{A{\displaystyle\prod\limits_k\left(1-\frac z{\mu_k}\right)}}{\displaystyle\pi^2z\prod\limits_k\left(1-\frac z{z_k(0)}\right)},\label{eq4.16}
\end{equation}
where number $A$ we find from condition \eqref{eq4.14}, and $\mu_k$ and $z_k(0)$ alternate. We will use the theorem by M. G. Krein \cite{9}.

\begin{theorem}[M. G. Krein]\label{t4.4}
In order that a real meromorphic function $f(z)$ map $\mathbb{C}_+\rightarrow\mathbb{C}_-$, it is necessary and sufficient that it is given by
\begin{equation}
f(z)=c\frac{z-a_0}{z-b_0}\prod\limits_{k\not=0}\left(1-\frac z{a_k}\right)\left(1-\frac z{b_k}\right)^{-1}\label{eq4.17}
\end{equation}
where $b_k<a_k<b_{k+1}$ ($k\in\mathbb{Z}_+$); $a_{-1}<0<b_1$; $c>0$.
\end{theorem}

Representation \eqref{eq4.17}, with $a_k<b_k<a_{k+1}$ ($\forall k\in\mathbb{Z}$) and $b_{-1}<0<a_1$, is also true for a real meromorphic function $f(z)$: $\mathbb{C}_+\rightarrow\mathbb{C}_-$. Along with the multiplicative representation \eqref{eq4.17}, there exists the additive decomposition \cite{9} for this class of meromorphic functions.

\begin{theorem}[N. G. Chebotarev]\label{t4.5}
Any real meromorphic function $f(z)$: $\mathbb{C}_+\rightarrow\mathbb{C}_+$ has the decomposition
\begin{equation}
f(z)=az+b-\frac{A_0}z+\sum\limits_{\omega_1}^{\omega_2}A_k\left(\frac1{b_k-z}-\frac1{b_k}\right)\label{eq4.18}
\end{equation}
where $-\infty\leq\omega_1\leq\omega_2\leq\infty$; $A_k\geq0$ ($\forall k\in\mathbb{Z}$); $a\geq0$, $b\in\mathbb{R}$, $b_k<b_{k+1}$ ($b_k\not=0$, $\forall k\in\mathbb{Z}$); the series
$$\sum\limits_{\omega_1}^{\omega_2}\frac{A_k}{b_k^2}<\infty$$
converges.
\end{theorem}

For real meromorphic functions $f(z)$: $\mathbb{C}_+\rightarrow\mathbb{C}_-$ there also exists decomposition \eqref{eq4.18} where $a\leq0$ and $A_k\leq0$ ($k\in\mathbb{Z}_+$).

Due to the Krein theorem, the function $F(z)$ \eqref{eq4.16} maps $\mathbb{C}_+\rightarrow\mathbb{C}_+$ (or $\mathbb{C}_+\rightarrow\mathbb{C}_-$). Using the Chebotarev theorem for $F(z)$, we obtain decomposition \eqref{eq4.18} where $\sign A_k=\sign A_0$. Relation \eqref{eq4.14} implies that $a=0$ and
$$b-\sum\frac{A_k}{z_k(0)}=1,$$
therefore
$$F(z)=1+\sum\limits_k\frac{A_k}{z_k(0)-z};\quad\sum\frac{|A_k|}{z_k(0)}<\infty.$$
Using \eqref{eq4.15}, we obtain
$$\sum\limits_k|A_k|<\infty.$$
Designating
$$\alpha=\sum_kA_k,\quad\|v_k\|^2=\frac{A_k}{\displaystyle\sum\limits_sA_s};$$
we have
$$F(z)=1+\alpha\sum\frac{\|v_k\|^2}{z_k(0)-z}.$$
It is left to set the function $v(x)$ such that $\|E_kv\|^2=\|v_k\|^2$ and to define the operator $L=L(\alpha,v)$ \eqref{eq1.9}. Its characteristic function $\Delta(\alpha,\lambda)$, in view of \eqref{eq2.14}, equals $\Delta(\alpha,\lambda)=F(\lambda^2)\Delta(0,\lambda)=\Delta(\lambda).$ $\blacksquare$

\renewcommand{\refname}{ \rm 

\newpage

Dr. Zolotarev V. A.

B. Verkin Institute for Low Temperature Physics and Engineering
of the National Academy of Sciences of Ukraine\\
47 Nauky Ave., Kharkiv, 61103, Ukraine

Department of Higher Mathematics and Informatics, V. N. Karazin Kharkov National University \\
4 Svobody Sq, Kharkov, 61077,  Ukraine

E-mail:vazolotarev@gmail.com

\end{Large}
\end{document}